%% file: Text-III.tex
\documentclass[11pt,twoside]{article}

\input{packages.tex}

\title{
  {\huge Complex Analysis of Real Functions \\[1.5ex]}
  III: Extended Fourier Theory }

\author{
  \Large Jorge L. deLyra\footnote{Email: delyra@latt.if.usp.br} \\
  Department of Mathematical Physics \\
  Physics Institute \\
  University of São Paulo }

\date{May 26, 2018}


\input{math-defs.tex}

\begin{document}\maketitle

\vspace{-1.65ex}
\begin{abstract}
  \noindent
  In the context of the complex-analytic structure within the unit disk
  centered at the origin of the complex plane, that was presented in a
  previous paper, we show that the complete Fourier theory of integrable
  real functions is contained within that structure, that is, within the
  structure of the space of inner analytic functions on the open unit
  disk. We then extend the Fourier theory beyond the realm of integrable
  real functions, to include for example singular Schwartz distributions,
  and possibly other objects.
\end{abstract}

\section{Introduction}\label{Sec01}

In a previous paper~\cite{CAoRFI} we introduced a certain complex-analytic
structure within the unit disk of the complex plane, and showed that one
can represent essentially all integrable real functions within that
structure. The construction leading to this result started with the use of
the Fourier coefficients $\alpha_{k}$ and $\beta_{k}$ of the integrable
real function $f(\theta)$, from which we defined a set of complex Taylor
coefficients $c_{k}$, thus leading to the corresponding inner analytic
function $w(z)$. It is therefore clearly apparent that there is a close
relation between that complex-analytic structure and the Fourier
theory~\cite{FSchurchill} of integrable real functions.

In this paper we will make that relation explicit by showing, in
Sections~\ref{Sec02}---\ref{Sec05}, that {\em all} the elements of the
Fourier theory of integrable real functions are contained within the
complex-analytic structure. What we mean by these elements is the set of
mathematical objects including the Fourier basis of functions, the Fourier
series, the scalar product for integrable real functions, the relations of
orthogonality and norm of the basis elements, and the completeness of the
Fourier basis, including its so-called completeness relation.

The fact that one can recover the real functions from their Fourier
coefficients almost everywhere, even when the corresponding Fourier series
are divergent, as we showed in~\cite{CAoRFI}, leads to a powerful and very
general summation rule for {\em all\/} Fourier series. Furthermore, we
will show in Section~\ref{Sec06} that the complex-analytic structure
allows us to extend the Fourier theory beyond the realm of integrable real
functions, to include the singular Schwartz distributions that we examined
in detail in another previous paper~\cite{CAoRFII}, as well as at least
some non-integrable real functions, and possibly other objects.

For ease of reference, we include here a one-page synopsis of the
complex-analytic structure introduced in~\cite{CAoRFI}. It consists of
certain elements within complex analysis~\cite{CVchurchill}, as well as of
their main properties.

\paragraph{Synopsis:} The Complex-Analytic Structure\\

\noindent
An {\em inner analytic function} $w(z)$ is simply a complex function which
is analytic within the open unit disk. An inner analytic function that has
the additional property that $w(0)=0$ is a {\em proper inner analytic
  function}. The {\em angular derivative} of an inner analytic function is
defined by

\noindent
\begin{equation}
  w^{\ldot}(z)
  =
  \ii
  z\,
  \frac{dw(z)}{dz}.
\end{equation}

\noindent
By construction we have that $w^{\ldot}(0)=0$, for all $w(z)$. The {\em
  angular primitive} of an inner analytic function is defined by

\begin{equation}
  w^{-1\ldot}(z)
  =
  -\ii
  \int_{0}^{z}dz'\,
  \frac{w(z')-w(0)}{z'}.
\end{equation}

\noindent
By construction we have that $w^{-1\ldot}(0)=0$, for all $w(z)$. In terms
of a system of polar coordinates $(\rho,\theta)$ on the complex plane,
these two analytic operations are equivalent to differentiation and
integration with respect to $\theta$, taken at constant $\rho$. These two
operations stay within the space of inner analytic functions, they also
stay within the space of proper inner analytic functions, and they are the
inverses of one another. Using these operations, and starting from any
proper inner analytic function $w^{0\ldot}(z)$, one constructs an infinite
{\em integral-differential chain} of proper inner analytic functions,

\begin{equation}
  \left\{
    \ldots,
    w^{-3\ldot}(z),
    w^{-2\ldot}(z),
    w^{-1\ldot}(z),
    w^{0\ldot}(z),
    w^{1\ldot}(z),
    w^{2\ldot}(z),
    w^{3\ldot}(z),
    \ldots\;
  \right\}.
\end{equation}

\noindent
Two different such integral-differential chains cannot ever intersect each
other. There is a {\em single} integral-differential chain of proper inner
analytic functions which is a constant chain, namely the null chain, in
which all members are the null function $w(z)\equiv 0$.

A general scheme for the classification of all possible singularities of
inner analytic functions is established. A singularity of an inner
analytic function $w(z)$ at a point $z_{1}$ on the unit circle is a {\em
  soft singularity} if the limit of $w(z)$ to that point exists and is
finite. Otherwise, it is a {\em hard singularity}. Angular integration
takes soft singularities to other soft singularities, and angular
differentiation takes hard singularities to other hard singularities.

Gradations of softness and hardness are then established. A hard
singularity that becomes a soft one by means of a single angular
integration is a {\em borderline hard} singularity, with degree of
hardness zero. The {\em degree of softness} of a soft singularity is the
number of angular differentiations that result in a borderline hard
singularity, and the {\em degree of hardness} of a hard singularity is the
number of angular integrations that result in a borderline hard
singularity. Singularities which are either soft or borderline hard are
integrable ones. Hard singularities which are not borderline hard are
non-integrable ones.

Given an integrable real function $f(\theta)$ on the unit circle, one can
construct from it a unique corresponding inner analytic function $w(z)$.
Real functions are obtained through the $\rho\to 1_{(-)}$ limit of the
real and imaginary parts of each such inner analytic function and, in
particular, the real function $f(\theta)$ is obtained from the real part
of $w(z)$ in this limit. The pair of real functions obtained from the real
and imaginary parts of one and the same inner analytic function are said
to be mutually Fourier-conjugate real functions.

Singularities of real functions can be classified in a way which is
analogous to the corresponding complex classification. Integrable real
functions are typically associated with inner analytic functions that have
singularities which are either soft or at most borderline hard. This ends
our synopsis.

\vspace{2.6ex}

\noindent
When we discuss real functions in this paper, some properties will be
globally assumed for these functions, just as was done in~\cite{CAoRFI}
and~\cite{CAoRFII}. These are rather weak conditions to be imposed on
these functions, that will be in force throughout this paper. It is to be
understood, without any need for further comment, that these conditions
are valid whenever real functions appear in the arguments. These weak
conditions certainly hold for any integrable real functions that are
obtained as restrictions of corresponding inner analytic functions to the
unit circle.

The most basic condition is that the real functions must be measurable in
the sense of Lebesgue, with the usual Lebesgue
measure~\cite{RARudin,RARoyden}. The second global condition we will
impose is that the functions have no removable singularities. The third
and last global condition is that the number of hard singularities on the
unit circle be finite, and hence that they be all isolated from one
another. There will be no limitation on the number of soft singularities.

The material contained in this paper is a development, reorganization and
extension of some of the material found, sometimes still in rather
rudimentary form, in the
papers~\cite{FTotCPI,FTotCPII,FTotCPIII,FTotCPIV,FTotCPV}.

\section{Fourier Series}\label{Sec02}

In~\cite{CAoRFI} we showed that, given any integrable real function
$f(\theta)$, one can construct a corresponding inner analytic function
$w(z)=u(\rho,\theta)+\ii v(\rho,\theta)$, from the real part of which
$f(\theta)$ can be recovered almost everywhere on the unit circle, through
the use of the $\rho\to 1_{(-)}$ limit, where $(\rho,\theta)$ are polar
coordinates on the complex plane. In that construction we started by
calculating the Fourier coefficients~\cite{FSchurchill} of the real
function, which is always possible given that the function is integrable,
using the usual integrals defining these coefficients,

\noindent
\begin{eqnarray}\label{EQFouCoe}
  \alpha_{0}
  & = &
  \frac{1}{\pi}
  \int_{-\pi}^{\pi}d\theta\,
  f(\theta),
  \nonumber\\
  \alpha_{k}
  & = &
  \frac{1}{\pi}
  \int_{-\pi}^{\pi}d\theta\,
  \cos(k\theta)f(\theta),
  \nonumber\\
  \beta_{k}
  & = &
  \frac{1}{\pi}
  \int_{-\pi}^{\pi}d\theta\,
  \sin(k\theta)f(\theta),
\end{eqnarray}

\noindent
for $k\in\{1,2,3,\ldots,\infty\}$. We then defined a set of complex Taylor
coefficients

\noindent
\begin{eqnarray}\label{EQTayCoe}
  c_{0}
  & = &
  \frac{1}{2}\,\alpha_{0},
  \nonumber\\
  c_{k}
  & = &
  \alpha_{k}
  -
  \ii\beta_{k},
\end{eqnarray}

\noindent
for $k\in\{1,2,3,\ldots,\infty\}$. Next we defined a complex variable $z$
associated to $\theta$, using the positive real variable $\rho$, by
$z=\rho\exp(\ii\theta)$. Using all these elements we then constructed the
power series

\begin{equation}\label{EQPowSer}
  S(z)
  =
  \sum_{k=0}^{\infty}
  c_{k}z^{k},
\end{equation}

\noindent
which we showed to be convergent to an inner analytic function $w(z)=S(z)$
within the open unit disk. This power series is therefore the Taylor
series of $w(z)$. We also proved that one recovers the real function
$f(\theta)$ almost everywhere on the unit circle from the $\rho\to
1_{(-)}$ limit of the real part $u(\rho,\theta)$ of $w(z)$. It is now very
easy to show that the Fourier series of an integrable real function
$f(\theta)$ is simply given by the real part of this Taylor series, when
restricted to the unit circle. Writing the series explicitly in terms of
the polar coordinates $(\rho,\theta)$ of the complex plane, we get

\noindent
\begin{eqnarray}\label{EQTayExp}
  w(z)
  & = &
  \frac{\alpha_{0}}{2}
  +
  \sum_{k=1}^{\infty}
  (\alpha_{k}-\ii\beta_{k})
  \rho^{k}
  \left[
    \cos(k\theta)+\ii\sin(k\theta)
  \right]
  \nonumber\\
  & = &
  \frac{\alpha_{0}}{2}
  +
  \sum_{k=1}^{\infty}
  \rho^{k}
  \left[
    \alpha_{k}\cos(k\theta)
    +
    \beta_{k}\sin(k\theta)
  \right]
  +
  \nonumber\\
  &   &
  \hspace{1.5em}
  +
  \ii
  \sum_{k=1}^{\infty}
  \rho^{k}
  \left[
    \alpha_{k}\sin(k\theta)
    -
    \beta_{k}\cos(k\theta)
  \right],
\end{eqnarray}

\noindent
where $w(z)=u(\rho,\theta)+\ii v(\rho,\theta)$. Taking now the $\rho\to
1_{(-)}$ limit we get

\noindent
\begin{eqnarray}
  u(1,\theta)
  +
  \ii
  v(1,\theta)
  & = &
  \frac{\alpha_{0}}{2}
  +
  \sum_{k=1}^{\infty}
  \left[
    \alpha_{k}\cos(k\theta)
    +
    \beta_{k}\sin(k\theta)
  \right]
  +
  \nonumber\\
  &   &
  \hspace{1.5em}
  +
  \ii
  \sum_{k=1}^{\infty}
  \left[
    \alpha_{k}\sin(k\theta)
    -
    \beta_{k}\cos(k\theta)
  \right].
\end{eqnarray}

\noindent
It follows, therefore, that the real part of $w(z)$ for $\rho=1$ is the
Fourier series of $f(\theta)$,

\begin{equation}
  u(1,\theta)
  =
  \frac{\alpha_{0}}{2}
  +
  \sum_{k=1}^{\infty}
  \left[
    \alpha_{k}\cos(k\theta)
    +
    \beta_{k}\sin(k\theta)
  \right],
\end{equation}

\noindent
and that the imaginary part of $w(z)$ for $\rho=1$ is the Fourier series
of the real function which is the Fourier conjugate of $f(\theta)$,

\begin{equation}
  v(1,\theta)
  =
  \sum_{k=1}^{\infty}
  \left[
    \alpha_{k}\sin(k\theta)
    -
    \beta_{k}\cos(k\theta)
  \right].
\end{equation}

\noindent
Here we see that, with respect to the Fourier series of $u(1,\theta)$, the
$k=0$ term is missing, all the other coefficients are the same, while the
$\cos(k\theta)$ were exchanged for $\sin(k\theta)$, and the
$\sin(k\theta)$ were exchanged for $-\cos(k\theta)$. In~\cite{CAoRFI} we
proved that $u(1,\theta)$ is equal to $f(\theta)$ almost everywhere,
irrespective of the convergence or lack of convergence of the Fourier
series, so that it now becomes clear that, when and where this
trigonometric series converges at all, it converges to the original
integrable real function,

\begin{equation}
  f(\theta)
  =
  \frac{\alpha_{0}}{2}
  +
  \sum_{k=1}^{\infty}
  \left[
    \alpha_{k}\cos(k\theta)
    +
    \beta_{k}\sin(k\theta)
  \right].
\end{equation}

\noindent
The convergence of this Fourier series can be characterized in terms of
the singularities of the inner analytic function $w(z)$ on the unit disk.
If there are no singularities of $w(z)$ on the unit circle, then the
maximum convergence disk of its Taylor series is larger than the unit
disk, and contains it. Therefore, in this case the Fourier series is
always convergent, as well as absolutely and uniformly convergent. On the
other hand, if there is at least one singularity of $w(z)$ on the unit
circle, then the unit disk is the maximum disk of convergence of the
Taylor series, and in this case the Fourier series may or may not be
convergent. In this case we see that, given any integrable real function,
the issue of the convergence of its Fourier series is thus identified
completely with the issue of the convergence of the Taylor series of the
corresponding inner analytic function, at the border of its maximum
convergence disk.

From the expansion in Equation~(\ref{EQTayExp}) we see that the recovery
of $f(\theta)$ from its Fourier coefficients via the inner analytic
function $w(z)$, as we discussed in~\cite{CAoRFI}, which works even when
the Fourier series diverges almost everywhere, is equivalent to taking the
$\rho\to 1_{(-)}$ limit of the following modified or {\em regulated}
Fourier series,

\noindent
\begin{eqnarray}\label{EQSummRule}
  f(\theta)
  & = &
  \lim_{\rho\to 1_{(-)}}
  \left\{
    \frac{\alpha_{0}}{2}
    +
    \sum_{k=1}^{\infty}
    \rho^{k}
    \left[
      \alpha_{k}\cos(k\theta)
      +
      \beta_{k}\sin(k\theta)
    \right]
  \right\}
  \nonumber\\
  & = &
  \frac{\alpha_{0}}{2}
  +
  \lim_{\rho\to 1_{(-)}}
  \sum_{k=1}^{\infty}
  \rho^{k}
  \left[
    \alpha_{k}\cos(k\theta)
    +
    \beta_{k}\sin(k\theta)
  \right],
\end{eqnarray}

\noindent
which of course is always convergent, so long as $\rho<1$, for all
integrable real functions $f(\theta)$, given that it is the real part of
the convergent Taylor series of $w(z)$. The limit indicated will exist
when and where $f(\theta)$ can be recovered from the real part of the
corresponding inner analytic function. This holds for all the points on
the unit circle where the inner analytic function $w(z)$ is either
analytic or has only soft singularities. This recipe constitutes,
therefore, a very general {\em summation rule} for Fourier series.

\section{Orthogonality Relations}\label{Sec03}

The Fourier series of an integrable real function can be understood as the
expansion of that real function in the Fourier basis of functions, which
consists of the set of functions

\begin{equation}\label{EQFouBas}
  \left\{
    \rule{0em}{2.5ex}
    1,
    \cos(k\theta),
    \sin(k\theta),
    k\in\{1,2,3,\ldots,\infty\}
  \right\}.
\end{equation}

\noindent
Let us now show that this is an orthogonal basis. Of course this can be
done using the standard form of the scalar product between two real
functions on the unit circle, by simply calculating a set of easy
integrals by elementary means. However, what we want to do here is to show
that both the form of the scalar product and the relations of
orthogonality and norm are contained within the structure of the inner
analytic functions, and can be derived from that structure. In fact, we
will show that these elements can be obtained from a particular set of
functions, the powers $z^{k}$, with $k\geq 0$, and their multiplicative
inverses $z^{-k}$. We start by noting that, if $C$ is any circle centered
at the origin, including the unit circle, then from the residues theorem
we have that

\begin{equation}
  \frac{1}{2\pi\ii}
  \oint_{C}dz\,
  z^{p-1}
  =
  \delta_{p,0},
\end{equation}

\noindent
where $p$ is an arbitrary integer, and where $\delta_{p,0}$ is the
Kronecker delta. This is so because the integral can be calculated by
residues, and a function which is a simple power, either positive or
negative, is its own Laurent series, which has only one term. Therefore,
its residue at $z=0$ is zero unless $p=0$, in which case it is equal to
one. Using this result for the case $p=k-k'$, where the integers $k$ and
$k'$ are in the set $\{0,1,2,3,\ldots,\infty\}$, we have

\begin{equation}
  \frac{1}{2\pi\ii}
  \oint_{C}dz\,
  z^{k-k'-1}
  =
  \delta_{k,k'},
\end{equation}

\noindent
while using the same expression for $p=k+k'$, with the limitation that
$k+k'>0$, which means that $k$ and $k'$ cannot both be zero, we have

\begin{equation}
  \frac{1}{2\pi\ii}
  \oint_{C}dz\,
  z^{k+k'-1}
  =
  0.
\end{equation}

\noindent
This is also a consequence of the Cauchy-Goursat theorem, since in this
case the integrand is analytic within the unit disk. Note that the power
$z^{k}$ with $k\geq 0$ is itself an inner analytic function. Writing these
two relations in terms of the integration variable $\theta$ we have

\noindent
\begin{eqnarray}
  \frac{1}{2\pi}\,
  \rho^{k-k'}
  \int_{-\pi}^{\pi}d\theta\,
  \e{\iii k\theta}
  \e{-\iii k'\theta}
  & = &
  \delta_{k,k'},
  \nonumber\\
  \frac{1}{2\pi}\,
  \rho^{k+k'}
  \int_{-\pi}^{\pi}d\theta\,
  \e{\iii k\theta}
  \e{\iii k'\theta}
  & = &
  0,
\end{eqnarray}

\noindent
since $z=\rho\exp(\ii\theta)$, where in the second equation we must have
$k+k'>0$. So long as $\rho\neq 0$ the powers of $\rho$ can be eliminated
from the second equation, and since the right-hand term of the first
equation is zero unless $k=k'$, they can also be eliminated from the first
equation, so that we have

\noindent
\begin{eqnarray}\label{EQComOrt}
  \frac{1}{2\pi}
  \int_{-\pi}^{\pi}d\theta\,
  \e{\iii k\theta}
  \e{-\iii k'\theta}
  & = &
  \delta_{k,k'},
  \nonumber\\
  \frac{1}{2\pi}
  \int_{-\pi}^{\pi}d\theta\,
  \e{\iii k\theta}
  \e{\iii k'\theta}
  & = &
  0,
\end{eqnarray}

\noindent
where in the second equation we must have $k+k'>0$. Note that this is
valid not only on the unit circle, but for all values of $\rho$ so long as
$\rho\neq 0$. Expanding the complex exponentials, with the use of the
Euler formula, and collecting real and imaginary parts, we have

\noindent
\begin{eqnarray}
  \frac{1}{\pi}
  \int_{-\pi}^{\pi}d\theta\,
  \left[
    \cos(k\theta)\cos(k'\theta)
    +
    \sin(k\theta)\sin(k'\theta)
  \right]
  \hspace{2em}
  &   &
  \nonumber\\
  +
  \ii\,
  \frac{1}{\pi}
  \int_{-\pi}^{\pi}d\theta\,
  \left[
    \sin(k\theta)\cos(k'\theta)
    -
    \cos(k\theta)\sin(k'\theta)
  \right]
  & = &
  2\delta_{k,k'},
  \nonumber\\
  \frac{1}{\pi}
  \int_{-\pi}^{\pi}d\theta\,
  \left[
    \cos(k\theta)\cos(k'\theta)
    -
    \sin(k\theta)\sin(k'\theta)
  \right]
  \hspace{2em}
  &   &
  \nonumber\\
  +
  \ii\,
  \frac{1}{\pi}
  \int_{-\pi}^{\pi}d\theta\,
  \left[
    \sin(k\theta)\cos(k'\theta)
    +
    \cos(k\theta)\sin(k'\theta)
  \right]
  & = &
  0,
\end{eqnarray}

\noindent
where in the second equation we must have $k+k'>0$. Since the right-hand
sides are real, we have the four real equations

\noindent
\begin{eqnarray}
  \frac{1}{\pi}
  \int_{-\pi}^{\pi}d\theta\,
  \left[
    \cos(k\theta)\cos(k'\theta)
    +
    \sin(k\theta)\sin(k'\theta)
  \right]
  & = &
  2\delta_{k,k'},
  \nonumber\\
  \frac{1}{\pi}
  \int_{-\pi}^{\pi}d\theta\,
  \left[
    \sin(k\theta)\cos(k'\theta)
    -
    \cos(k\theta)\sin(k'\theta)
  \right]
  & = &
  0,
  \nonumber\\
  \frac{1}{\pi}
  \int_{-\pi}^{\pi}d\theta\,
  \left[
    \cos(k\theta)\cos(k'\theta)
    -
    \sin(k\theta)\sin(k'\theta)
  \right]
  & = &
  0,
  \nonumber\\
  \frac{1}{\pi}
  \int_{-\pi}^{\pi}d\theta\,
  \left[
    \sin(k\theta)\cos(k'\theta)
    +
    \cos(k\theta)\sin(k'\theta)
  \right]
  & = &
  0,
\end{eqnarray}

\noindent
where we must have $k+k'>0$ in the last two equations. In the case
$k+k'=0$, which implies that $k=0$ and $k'=0$, we obtain from the first
equation the identity

\begin{equation}\label{EQCaseZero}
  \frac{1}{\pi}
  \int_{-\pi}^{\pi}d\theta\,
  \cos(0)\cos(0)
  =
  2,
\end{equation}

\noindent
which is a part of the relations of orthogonality and norm of the Fourier
basis, namely the one giving the squared norm of the constant function
which is equal to one for all $\theta$. The second equation is just a
trivial identity when we have $k=0$ and $k'=0$, which we may therefore
ignore. We may now assume that we have $k+k'>0$ for all the four
equations. Adding and subtracting the first and third equations we get

\noindent
\begin{eqnarray}
  \frac{1}{\pi}
  \int_{-\pi}^{\pi}d\theta\,
  \cos(k\theta)\cos(k'\theta)
  & = &
  \delta_{k,k'},
  \nonumber\\
  \frac{1}{\pi}
  \int_{-\pi}^{\pi}d\theta\,
  \sin(k\theta)\sin(k'\theta)
  & = &
  \delta_{k,k'},
\end{eqnarray}

\noindent
for $k+k'>0$, while adding and subtracting the other two equations we get

\noindent
\begin{eqnarray}
  \frac{1}{\pi}
  \int_{-\pi}^{\pi}d\theta\,
  \sin(k\theta)\cos(k'\theta)
  & = &
  0,
  \nonumber\\
  \frac{1}{\pi}
  \int_{-\pi}^{\pi}d\theta\,
  \cos(k\theta)\sin(k'\theta)
  & = &
  0,
\end{eqnarray}

\noindent
for $k+k'>0$, which are just two copies of the same relation. We have
therefore the complete set of orthogonality relations, which also includes
those relations giving the norms of the basis functions,

\noindent
\begin{eqnarray}
  \frac{1}{\pi}
  \int_{-\pi}^{\pi}d\theta\,
  \cos(k\theta)\cos(k'\theta)
  & = &
  \delta_{k,k'},
  \nonumber\\
  \frac{1}{\pi}
  \int_{-\pi}^{\pi}d\theta\,
  \sin(k\theta)\sin(k'\theta)
  & = &
  \delta_{k,k'},
  \nonumber\\
  \frac{1}{\pi}
  \int_{-\pi}^{\pi}d\theta\,
  \sin(k\theta)\cos(k'\theta)
  & = &
  0,
\end{eqnarray}

\noindent
where $k+k'>0$, which includes all the relevant cases, that is, all the
relevant pairs of elements of the basis in Equation~(\ref{EQFouBas}),
except for the single case for $k=0$ and $k'=0$, which we examined
separately before, leading to Equation~(\ref{EQCaseZero}). Note that this
derivation included the determination of the form of the scalar product
for the basis elements. Given two integrable real functions $f(\theta)$
and $g(\theta)$, their scalar product is given by

\begin{equation}\label{EQScalProd}
  (f|g)
  =
  \int_{-\pi}^{\pi}d\theta\,
  f(\theta)g(\theta),
\end{equation}

\noindent
which induces a positive-definite norm in the space of all integrable real
functions defined on the periodic interval, which is thus seen to
constitute a Hilbert space. We may therefore conclude that the whole
structure of orthogonality and norm of the Fourier basis is contained in
the structure of the inner analytic function within the unit disk of the
complex plane.

Note that, since all possible inner analytic functions are given by
convergent power series within the open unit disk, and since these power
series can be understood as infinite linear combinations of the particular
set of inner analytic functions given by the non-negative powers
$\left\{\rule{0em}{2ex}z^{k},k\in\{0,1,2,3,\ldots,\infty\}\right\}$, we
may think that this set of functions forms a {\em basis} of the space of
inner analytic functions, which we may call the {\em Taylor basis}. Since
the orthogonality of the Fourier basis was obtained above from the
properties of this set of non-negative powers, it becomes clear that the
orthogonality of the Fourier basis is a consequence of similar properties
that must hold for the Taylor basis. In fact, it is possible to define a
complex scalar product within the space of inner analytic functions,
according to which this Taylor basis is orthogonal. Since this constitutes
a considerable detour from our main line of reasoning here, it will be
presented as an appendix. As one can see in Appendix~\ref{App01}, this
complex scalar product induces in the space of inner analytic functions a
positive-definite norm. As was observed in~\cite{CAoRFI}, this space forms
a vector space over the field of complex numbers, and we thus see that it
constitutes in fact a complex Hilbert space.

\section{Completeness Relation}\label{Sec04}

Let us now prove the completeness of the Fourier basis. In this context
the concept of completeness is that of a basis within a vector space. We
will first give a simple and direct proof of completeness, which is
however subject to a slight limitation regarding the vector space for
which the basis is shown to be complete, using the analytic structure
within the open unit disk, and later establish the relation of the concept
of completeness with the so-called completeness relation. The proof of
completeness using the completeness relation is not subject to any such
limitation.

\vspace{2.6ex}

\noindent
In this section we will prove the following completeness theorem.

\begin{theorem}\Colon\label{Theo01}
  The basis of real functions
  $\left\{\rule{0em}{2ex}1,\cos(k\theta),\sin(k\theta),
    k\in\{1,2,3,\ldots,\infty\}\right\}$, is complete to represent the
  space of all integrable real functions defined on the unit circle.
\end{theorem}

\noindent
The proof consists of establishing that, given an arbitrary integrable
real function $\psi(\theta)$ on the unit circle, which is orthogonal to
all the elements of the Fourier basis, according to the scalar product
defined in Equation~(\ref{EQScalProd}), it then follows that
$\psi(\theta)$ must be zero almost everywhere. Note that the orthogonality
to the elements of the basis means that $\psi(\theta)$ is such that all
its Fourier coefficients, as defined in Equation~(\ref{EQFouCoe}), are
zero.

\begin{proof}\Colon
\end{proof}

\noindent
Let $\psi(\theta)$ be a real function on the unit circle which can be
obtained as the $\rho\to 1_{(-)}$ limit of an inner analytic function. We
assume that it is orthogonal to all the elements of the basis, so that all
its Fourier coefficients are zero, that is, we assume that for this
function we have $\alpha_{0}=0$, $\alpha_{k}=0$ and $\beta_{k}=0$, for all
$k\in\{1,2,3,\ldots,\infty\}$. Since we thus have all the Fourier
coefficients of $\psi(\theta)$, we may use the construction presented
in~\cite{CAoRFI} in order to determine the corresponding inner analytic
function. However, since all the Fourier coefficients are zero, it follows
at once from the step of that construction given in
Equation~(\ref{EQTayCoe}) that for $\psi(\theta)$ the complex coefficients
$c_{k}$ are zero for all $k$. Therefore, the power series $S(z)$
constructed in the next step of the process, given in
Equation~(\ref{EQPowSer}), is identically zero and thus converges
trivially to the identically zero complex function $w_{\psi}(z)\equiv 0$
on the whole complex plane.

The analyticity region of $w_{\psi}(z)$ includes the unit circle, and
therefore the series converges to zero there. Since on the one hand the
series converges to zero, and on the other hand we know that for $\rho=1$
it necessarily converges to the restriction of $w_{\psi}(z)$ to the unit
circle, it follows that the restriction, including both real and imaginary
parts, must be zero everywhere on the unit circle. Therefore it follows
that $\psi(\theta)$ and the identically zero real function coincide
everywhere on the unit circle, and therefore we conclude that
$\psi(\theta)=0$ everywhere on that circle. This establishes that the
Fourier basis is complete for the space of all integrable real functions
defined on the periodic interval, which can be obtained as the $\rho\to
1_{(-)}$ limits of inner analytic functions. This completes the first
version of the proof of Theorem~\ref{Theo01}, which is valid for the
vector space of real functions just described.

\vspace{2.6ex}

\noindent
Note that, since all possible inner analytic functions are given by
convergent power series within the open unit disk, and since these power
series can be understood as expansions of those inner analytic functions
in the Taylor basis of functions given by the non-negative powers
$\left\{\rule{0em}{2ex}z^{k},k\in\{0,1,2,3,\ldots,\infty\}\right\}$, we
may say that this Taylor basis is complete for the space of all inner
analytic functions. Since the proof of the completeness of the Fourier
basis given above was obtained from the complex-analytic structure within
the open unit disk, it becomes clear that the completeness of the Fourier
basis on the unit circle is a consequence of the completeness of the
Taylor basis within the open unit disk. This adds to the relationship
between the Fourier basis on the unit circle and the Taylor basis on the
unit disk, which was first established during the discussion involving the
orthogonality of the Fourier basis, in Section~\ref{Sec03}. In addition to
all this, within the spaces generated by either basis one may define
scalar products that induce positive-definite norms, thus making them both
Hilbert spaces, as is discussed in Appendix~\ref{App01}.

Let us now turn to the usual completeness relation. Let us first write it
down and then exhibit its usefulness. The relation can be understood as
the expression, as a Fourier series, of the Dirac delta ``function''
defined with respect to a point given by the angle $\theta_{1}$ on the
unit circle, which we examined in great detail in~\cite{CAoRFII}, and
which we denote by $\delta(\theta-\theta_{1})$. As we have shown
in~\cite{CAoRFII}, using the usual rules for the manipulation of the delta
``function'', one finds that the corresponding Fourier coefficients are
given by

\noindent
\begin{eqnarray}
  \alpha_{0}
  & = &
  \frac{1}{\pi}
  \int_{-\pi}^{\pi}d\theta\,
  \delta(\theta-\theta_{1})
  \nonumber\\
  & = &
  \frac{1}{\pi},
  \nonumber\\
  \alpha_{k}
  & = &
  \frac{1}{\pi}
  \int_{-\pi}^{\pi}d\theta\,
  \cos(k\theta)\delta(\theta-\theta_{1})
  \nonumber\\
  & = &
  \frac{1}{\pi}\,
  \cos(k\theta_{1}),
  \nonumber\\
  \beta_{k}
  & = &
  \frac{1}{\pi}
  \int_{-\pi}^{\pi}d\theta\,
  \sin(k\theta)\delta(\theta-\theta_{1})
  \nonumber\\
  & = &
  \frac{1}{\pi}\,
  \sin(k\theta_{1}),
\end{eqnarray}

\noindent
for $k\in\{1,2,3,\ldots,\infty\}$, so that the completeness relation is
given by the Fourier expansion, that turns out to be a bi-linear form on
the elements of the Fourier basis,

\begin{equation}\label{EQExpanDel}
  \delta(\theta-\theta_{1})
  =
  \frac{1}{2\pi}
  +
  \frac{1}{\pi}
  \sum_{k=1}^{\infty}
  \left[
    \cos(k\theta_{1})
    \cos(k\theta)
    +
    \sin(k\theta_{1})
    \sin(k\theta)
  \right],
\end{equation}

\noindent
which is manifestly divergent, but which can be made to converge for all
values of $\theta$, so that we may recover the delta ``function'' almost
everywhere, in fact everywhere but at $\theta_{1}$, through the use of the
summation rule given in Equation~(\ref{EQSummRule}),

\begin{equation}\label{EQRegExpDel}
  \delta(\theta-\theta_{1})
  =
  \frac{1}{2\pi}
  +
  \frac{1}{\pi}
  \lim_{\rho\to 1_{(-)}}
  \sum_{k=1}^{\infty}
  \rho^{k}
  \left[
    \cos(k\theta_{1})
    \cos(k\theta)
    +
    \sin(k\theta_{1})
    \sin(k\theta)
  \right].
\end{equation}

\noindent
This is equivalent to the definition of the delta ``function'' as the
$\rho\to 1_{(-)}$ limit of the real part of the inner analytic function
given by

\begin{equation}
  w_{\delta}(z,z_{1})
  =
  \frac{1}{2\pi}
  -
  \frac{1}{\pi}\,
  \frac{z}{z-z_{1}},
\end{equation}

\noindent
as was discussed in detail in~\cite{CAoRFII}. One can use the expansion in
Equation~(\ref{EQExpanDel}), possibly regulated as in
Equation~(\ref{EQRegExpDel}), to prove the completeness of the basis,
while operating strictly in terms of real objects on or near the unit
circle. Here is how this can be done.

\newpage

\begin{proof}\Colon
\end{proof}

\noindent
If we assume that an arbitrary integrable real function $\psi(\theta)$ on
the unit circle is given, which is such that its scalar products with all
the elements of the basis are zero, then we have the infinite set of
equations

\noindent
\begin{eqnarray}
  \int_{-\pi}^{\pi}d\theta\,
  \psi(\theta)
  & = &
  0,
  \nonumber\\
  \int_{-\pi}^{\pi}d\theta\,
  \cos(k\theta)\psi(\theta)
  & = &
  0,
  \nonumber\\
  \int_{-\pi}^{\pi}d\theta\,
  \sin(k\theta)\psi(\theta)
  & = &
  0,
\end{eqnarray}

\noindent
for all $k\in\{1,2,3,\ldots,\infty\}$. We may therefore construct an
infinite linear combination of all these equations, with the coefficients
carefully chosen as shown below, involving an arbitrary parameter
$\theta_{1}$ in the interval $[-\pi,\pi]$ and an auxiliary strictly
positive real variable $\rho<1$, where the right-hand side is still zero,

\noindent
\begin{eqnarray*}
  \left[
    \frac{1}{2\pi}
  \right]
  \int_{-\pi}^{\pi}d\theta\,
  \psi(\theta)
  +
  \sum_{k=1}^{\infty}
  \left[
    \rho^{k}
    \frac{1}{\pi}
    \cos(k\theta_{1})
  \right]
  \int_{-\pi}^{\pi}d\theta\,
  \cos(k\theta)\psi(\theta)
  +
  &   &
  \\
  +
  \sum_{k=1}^{\infty}
  \left[
    \rho^{k}
    \frac{1}{\pi}
    \sin(k\theta_{1})
  \right]
  \int_{-\pi}^{\pi}d\theta\,
  \sin(k\theta)\psi(\theta)
  \hspace{1em}
  & = &
  0
  \;\;\;\Rightarrow
\end{eqnarray*}
\begin{equation}
  \int_{-\pi}^{\pi}d\theta\,
  \left\{
    \frac{1}{2\pi}
    +
    \frac{1}{\pi}
    \sum_{k=1}^{\infty}
    \rho^{k}
    \left[
      \cos(k\theta_{1})
      \cos(k\theta)
      +
      \sin(k\theta_{1})
      \sin(k\theta)
    \right]
  \right\}
  \psi(\theta)
  \;\;=\;\;
  0.
\end{equation}

\noindent
Since the expression within curly brackets in this last integral is now
seen to be the regulated expansion of $\delta(\theta-\theta_{1})$ in the
Fourier basis, shown in Equation~(\ref{EQRegExpDel}), we may therefore
take the $\rho\to 1_{(-)}$ limit and write that

\begin{equation}
  \int_{-\pi}^{\pi}d\theta\,
  \delta(\theta-\theta_{1})
  \psi(\theta)
  =
  0.
\end{equation}

\noindent
Finally, using the rules of manipulation of the delta ``function'', when
and where $\psi(\theta)$ is continuous, which it therefore must be almost
everywhere, we have

\begin{equation}
  \psi(\theta_{1})
  =
  0.
\end{equation}

\noindent
Since $\theta_{1}$ is an arbitrary value of $\theta$, we conclude that
$\psi(\theta)$ is zero everywhere. This completes the second version of
the proof of Theorem~\ref{Theo01}, which is valid for the vector space of
all integrable real functions defined on the unit circle, regardless of
whether or not they can be obtained from an inner analytic function.

\vspace{2.6ex}

\noindent
Note that, in a sense, this method of proof of the completeness of the
Fourier basis is a little more limited than the direct proof using the
analytic structure within the open unit disk, because we must assume
during the argument that $\psi(\theta)$ is continuous almost everywhere.
However, since this hypothesis does get confirmed a posteriori by the
result obtained, this is not a true limitation.

On the other hand, this second proof is {\em less} limited than the first
one because in this case the vector space of functions for which one shows
that the basis is complete is the space of integrable real functions
without removable singularities defined on the interval $[-\pi,\pi]$, with
{\em no reference} to whether or not these functions can be obtained as
the $\rho\to 1_{(-)}$ limits of inner analytic functions.

In fact, by establishing the completeness of the Fourier basis without any
recourse to the $\rho\to 1_{(-)}$ limit for the real functions, as a
corollary of this second proof we have shown that there is no integrable
real function on the unit circle, other that the identically zero real
function, which corresponds to the identically zero inner analytic
function. As a consequence of this, there is no integrable real function
defined on the unit circle that cannot be represented by an unique inner
analytic function.

\section{Notes on the Convergence Problem}\label{Sec05}

In this paper we have made the deliberate choice of not discussing the
question of the convergence of Fourier series in any amount of detail,
that is, we have not discussed any of the many existing so-called Fourier
theorems. The reason for this is that we believe that this would
constitute a rather long and complex discussion, best left for a separate
paper. Instead, we have focused our attention on the summation rule given
in Equation~(\ref{EQSummRule}), according to which {\em all\/} Fourier
series of integrable real functions, without any further restrictions, can
be added up in such a way that one is able to recover the functions from
their Fourier coefficients, even if the Fourier series themselves diverge.
However, we may make a few comments about the issue of convergence,
without going too far afield in that subject, in order to exhibit the
relation between our complex analytic structure and the convergence
problem.

First of all, let us recall that, as was shown in~\cite{CAoRFI}, the real
function $f(\theta)$ is equal almost everywhere to the real part of the
corresponding inner analytic function $w(z)$, taken in the $\rho\to
1_{(-)}$ limit, and also that, as we have shown in Section~\ref{Sec02} of
this paper, the Fourier series of $f(\theta)$ is given by the real part of
the Taylor series $S(z)$ of $w(z)$ in that same limit. Therefore, it is
clearly apparent that, as was already noted in Section~\ref{Sec02}, the
problem of the convergence of real Fourier series is completely identified
with the problem of the convergence of the corresponding complex power
series on the unit circle, including the cases in which it is the rim of
their maximum disks of convergence. Whatever is established for one type
of series is also valid for the other. As was also noted in
Section~\ref{Sec02}, the convergence properties on the unit circle will
depend on the existence and nature of the singularities of $w(z)$ on that
circle.

One way to discuss the issue of convergence is to observe that the
summation rule given in Equation~(\ref{EQSummRule}) involves two limits,
one being the series summation limit and the other being the $\rho\to
1_{(-)}$ limit from the interior of the unit disk to the unit circle. What
has been shown so far in this series of papers is that if one takes the
series summation limit first, and only after that the $\rho\to 1_{(-)}$
limit, then it is always possible to recover the real function from its
Fourier coefficients. It is therefore immediately apparent that the
statement that the Fourier series converges over the unit circle is
equivalent to the statement that the order of these two limits can be
inverted. In fact, by first taking the $\rho\to 1_{(-)}$ limit one obtains
the usual Fourier series over the unit circle, and if one is then able to
take the series summation limit, then that series converges to the
corresponding real function.

The general problem of deciding under what conditions the order of the two
limits can be inverted is not a simple one. However, it is not too
difficult to use our analytic structure to write the partial sums of the
Fourier series in terms of real integrals which are similar to the
Dirichlet integrals usually involved in some of the Fourier theorems.
This can then be used as the starting point for further discussions of the
convergence problem, including in particular discussions establishing the
connection of the analytic structure with specific Fourier theorems. In
order to do this, let $f(\theta)$ be an integrable real function on
$[-\pi,\pi]$ and let the real numbers $\alpha_{0}$, $\alpha_{k}$ and
$\beta_{k}$, for $k\in\{1,2,3,\ldots,\infty\}$, be its Fourier
coefficients. We may define the complex coefficients $c_{0}$ and $c_{k}$
shown in Equation~(\ref{EQTayCoe}), and thus construct the corresponding
inner analytic function $w(z)$ within the open unit disk, using the power
series $S(z)$ given in Equation~(\ref{EQPowSer}), which, as was shown
in~\cite{CAoRFI}, always converges for $|z|<1$. The partial sums of the
first $N$ terms of this series are given by

\begin{equation}
  S_{N}(z)
  =
  \sum_{k=0}^{N-1}
  c_{k}z^{k},
\end{equation}

\noindent
a complex sequence which, for $|z|<1$, we already know to converge to
$w(z)$ in the $N\to\infty$ limit. Note however that, since $S_{N}(z)$ is
in fact an analytic function over the whole complex plane, this expression
itself can be consistently considered for all $z$, and in particular for
$z$ on the unit circle, where $|z|=1$. Note also that the function $w(z)$
may have singularities on the unit circle, but that these must be
integrable ones, at least along that circle. In addition to this, the
complex coefficients $c_{k}$ may be written as integrals involving $w(z)$,
with the use of the Cauchy integral formulas,

\begin{equation}
  c_{k}
  =
  \frac{1}{2\pi\ii}
  \oint_{C}dz\,
  \frac{w(z)}{z^{k+1}},
\end{equation}

\noindent
where $C$ can be taken as a circle centered at the origin, with radius
$\rho\leq 1$. The reason why we may include the case $\rho=1$ here is
that, as was shown in~\cite{CAoRFI}, as a function of $\rho$ the
expression above for $c_{k}$ is not only constant within the open unit
disk, but also continuous from within at the unit circle. In this way the
coefficients $c_{k}$ may be written back in terms of the inner analytic
function $w(z)$. If we substitute this expression for $c_{k}$ back in the
partial sums of the series we get

\noindent
\begin{eqnarray}
  S_{N}(z)
  & = &
  \sum_{k=0}^{N-1}
  z^{k}\,
  \frac{1}{2\pi\ii}
  \oint_{C}dz_{1}\,
  \frac{w(z_{1})}{z_{1}^{k+1}}
  \nonumber\\
  & = &
  \frac{1}{2\pi\ii}
  \oint_{C}dz_{1}\,
  \frac{w(z_{1})}{z_{1}}
  \sum_{k=0}^{N-1}
  \left(
    \frac{z}{z_{1}}
  \right)^{k},
\end{eqnarray}

\noindent
where we must have $|z_{1}|\leq 1$. The sum is now a finite geometric
progression, so that we have

\noindent
\begin{eqnarray}\label{EQParSum}
  S_{N}(z)
  & = &
  \frac{1}{2\pi\ii}
  \oint_{C}dz_{1}\,
  \frac{w(z_{1})}{z_{1}}\,
  \frac
  {1-(z/z_{1})^{N}}
  {1-(z/z_{1})}
  \nonumber\\
  & = &
  \frac{1}{2\pi\ii}
  \oint_{C}dz_{1}\,
  \frac{w(z_{1})}{z_{1}-z}
  -
  \frac{z^{N}}{2\pi\ii}
  \oint_{C}dz_{1}\,
  \frac{w(z_{1})}{z_{1}^{N}(z_{1}-z)}.
\end{eqnarray}

\noindent
A careful discussion of this formula is now in order. There are two
relevant cases to consider. In the first case we see that, if we have that
$|z_{1}|>|z|$, then in the first term above we obtain the expression of
the Cauchy integral formula for $w(z)$, which then allows us to write an
explicit expression for the remainder of the complex power series after
one adds up its first $N$ terms,

\noindent
\begin{eqnarray}
  R_{N}(z)
  & = &
  w(z)
  -
  S_{N}(z)
  \nonumber\\
  & = &
  \frac{z^{N}}{2\pi\ii}
  \oint_{C}dz_{1}\,
  \frac{w(z_{1})}{z_{1}^{N}(z_{1}-z)},
\end{eqnarray}

\noindent
where $|z|<|z_{1}|\leq 1$. This expression of the remainder in closed
form, an expression which, as one can easily show, goes to zero in the
$N\to\infty$ limit, is what makes it easy to discuss the convergence of
complex power series. However, this expression does {\em not} give us an
equivalent expression for the remainder of the Fourier series, because
this would require us to make $|z|=|z_{1}|=1$, which is not allowed by the
strict inequality $|z|<|z_{1}|$, a restriction which is due to the use of
the Cauchy integral formulas. In the second case we observe that, if we
have that $|z_{1}|<|z|$, then the first term in Equation~(\ref{EQParSum})
is simply zero, and therefore we get a modified expression for the partial
sums of the series,

\begin{equation}
  S_{N}(z)
  =
  -\,
  \frac{z^{N}}{2\pi\ii}
  \oint_{C}dz_{1}\,
  \frac{w(z_{1})}{z_{1}^{N}(z_{1}-z)},
\end{equation}

\noindent
where $|z|>|z_{1}|$. Note that in this case we are unable to write an
explicit expression in closed form for the remainder of the series, a fact
which seems to be related to the remarkable difficulty in finding a
necessary and sufficient condition for the convergence of Fourier series.
Since $z$ may have any complex value in this expression, we may now make
$z=\rho\exp(\ii\theta)$ with $\rho=1$, as well as
$z_{1}=\rho_{1}\exp(\ii\theta_{1})$, and thus write the integral
explicitly in terms of the variable $\theta_{1}$ on the circle of radius
$\rho_{1}$,

\noindent
\begin{eqnarray}
  S_{N}(1,\theta)
  & = &
  -\,
  \frac{\e{\iii N\theta}}{2\pi\ii}
  \int_{-\pi}^{\pi}d\theta_{1}\,
  \ii\rho_{1}\e{\iii\theta_{1}}\,
  \frac
  {
    w(\rho_{1},\theta_{1})
  }
  {
    \rho_{1}^{N}\e{\iii N\theta_{1}}
    \left(
      \rho_{1}\e{\iii\theta_{1}}-\e{\iii\theta}
    \right)
  }
  \nonumber\\
  & = &
  -\,
  \frac{1}{2\pi\rho_{1}^{N-1}}
  \int_{-\pi}^{\pi}d\theta_{1}\,
  \e{\iii N(\theta-\theta_{1})}\,
  \frac
  {
    w(\rho_{1},\theta_{1})
  }
  {
    \rho_{1}-\e{\iii(\theta-\theta_{1})}
  }.
\end{eqnarray}

\noindent
Making $\Delta\theta=\theta_{1}-\theta$ we have

\begin{equation}
  S_{N}(1,\theta)
  =
  -\,
  \frac{1}{2\pi\rho_{1}^{N-1}}
  \int_{-\pi}^{\pi}d\theta_{1}\,
  \e{-\iii N\Delta\theta}\,
  \frac
  {
    w(\rho_{1},\theta_{1})
  }
  {
    \rho_{1}-\e{-\iii\Delta\theta}
  }.
\end{equation}

\noindent
In order to be able to write explicitly the real and imaginary parts of
the partial sums, we must now rationalize this expression,

\noindent
\begin{eqnarray}
  S_{N}(1,\theta)
  & = &
  -\,
  \frac{1}{2\pi\rho_{1}^{N-1}}
  \int_{-\pi}^{\pi}d\theta_{1}\,
  \e{-\iii N\Delta\theta}\,
  \frac
  {
    w(\rho_{1},\theta_{1})
    \left(
      \rho_{1}-\e{\iii\Delta\theta}
    \right)
  }
  {
    \left(
      \rho_{1}-\e{-\iii\Delta\theta}
    \right)
    \left(
      \rho_{1}-\e{\iii\Delta\theta}
    \right)
  }
  \nonumber\\
  & = &
  \frac{1}{2\pi\rho_{1}^{N-1}}
  \int_{-\pi}^{\pi}d\theta_{1}\,
  w(\rho_{1},\theta_{1})\,
  \e{-\iii (N-1/2)\Delta\theta}
  \times
  \nonumber\\
  &   &
  \hspace{7em}
  \times\,
  \frac
  {
    \e{\iii\Delta\theta/2}
    -
    \rho_{1}
    \e{-\iii\Delta\theta/2}
  }
  {
    1+\rho_{1}^{2}-2\rho_{1}\cos(\Delta\theta)
  }.
\end{eqnarray}

\noindent
The expression can be somewhat simplified if we write most things in terms
of $\Delta\theta/2$, as well as in terms of $N_{1}=N-1/2$,

\noindent
\begin{eqnarray}
  S_{N}(1,\theta)
  & = &
  \frac{1}{2\pi\rho_{1}^{N-1}}
  \int_{-\pi}^{\pi}d\theta_{1}\,
  \left[
    \rule{0em}{2.5ex}
    u(\rho_{1},\theta_{1})
    +
    \ii
    v(\rho_{1},\theta_{1})
  \right]
  \times
  \nonumber\\
  &   &
  \hspace{7em}
  \times
  \left[
    \rule{0em}{2.5ex}
    \cos(N_{1}\Delta\theta)
    -
    \ii
    \sin(N_{1}\Delta\theta)
  \right]
  \times
  \\
  &   &
  \hspace{7em}
  \times\,
  \frac
  {
    (1-\rho_{1})\cos(\Delta\theta/2)
    +
    \ii
    (1+\rho_{1})\sin(\Delta\theta/2)
  }
  {
    (1-\rho_{1})^{2}+4\rho_{1}\sin^{2}(\Delta\theta/2)
  }.
  \nonumber
\end{eqnarray}

\noindent
In this context, a Fourier theorem is one which states sufficient
conditions on $f(\theta)$ under which it follows that the real part of the
corresponding sequence of partial sums $S_{N}(1,\theta)$ converges in the
$N\to\infty$ limit, after one takes the $\rho_{1}\to 1$ limit, so that the
integral is written over the unit circle. In any circumstances in which
one managed to calculate these integrals explicitly in terms of
$\rho_{1}$, for $\rho_{1}<1$, one would then be able to consider taking
the $\rho_{1}\to 1$ limit of the resulting expression. However, despite
the facts that $f(\theta)=u(1,\theta_{1})$ and that
$g(\theta)=v(1,\theta_{1})$, almost everywhere over the unit circle, as
well as the fact that these are integrable real functions, we cannot
simply take the $\rho_{1}\to 1$ limit of this expression as it stands,
because it was derived under the hypothesis that $|z|>|z_{1}|$, and
therefore that $\rho>\rho_{1}$, which at this point implies the strict
inequality $1>\rho_{1}$. We may however put $\rho_{1}=1$ in the integrand
simply in order to simplify the integrals, so as to exhibit their
structure more clearly. If one does that one obtains

\begin{equation}
  \int_{-\pi}^{\pi}d\theta_{1}\,
  \left[
    \rule{0em}{2.5ex}
    f(\theta_{1})
    +
    \ii
    g(\theta_{1})
  \right]\,
  \frac
  {
    \sin\!\left[\rule{0em}{2ex}(N-1/2)\Delta\theta\right]
    +
    \ii
    \cos\!\left[\rule{0em}{2ex}(N-1/2)\Delta\theta\right]
  }
  {
    \sin(\Delta\theta/2)
  },
\end{equation}

\noindent
which clearly reduces to Dirichlet integrals and other similar integrals.
A more complete discussion of the issue of convergence would require
considerable development of the ideas and structures involved in these
arguments. It is currently not entirely clear how useful the analytic
structure within the open unit disk can be in regards to proving known
Fourier theorems or discovering new ones.

\section{Extension of the Theory}\label{Sec06}

Up to this point we have been examining only the Fourier theory of
integrable real functions. In addition to this, a small extension of the
theory has already been considered when we wrote the Fourier expansion of
the Dirac delta ``function'' in Equations~(\ref{EQExpanDel})
and~(\ref{EQRegExpDel}) of Section~\ref{Sec04}, with the help of the
summation rule given in Equation~(\ref{EQSummRule}). This ``function'' has
in common with the integrable real functions the fact that its Fourier
coefficients $\alpha_{k}$ and $\beta_{k}$ are limited when we take the
limit $k\to\infty$. The same is true for the corresponding complex Taylor
coefficients $c_{k}$ in either case. However, the correspondence between
real Fourier coefficients and complex Taylor coefficients given by the
relations in Equation~(\ref{EQTayCoe}) can be generalized, independently
of any concerns about the behavior of these coefficients when
$k\to\infty$, and independently of any concerns about the convergence of
the corresponding series.

We will now discuss the extension of the Fourier theory beyond the realm
of integrable real functions. One way to look at this, which is probably
the most general possible way, is to simply consider the set of {\em all}
inner analytic functions. Given {\em any} inner analytic function $w(z)$
and its complex Taylor series around the origin, which is therefore
convergent within the open unit disk, and irrespective of whether or not
$w(z)$ corresponds to an integrable real function, one can define a
corresponding real Fourier series on the unit circle. In all such cases
the issues of convergence of the resulting Fourier series are then
completely identified with the corresponding issues for the Taylor series
restricted to the unit circle, which is often the border of its maximum
disk of convergence. Important examples which are not related to
integrable real functions are the cases of the Dirac delta ``function''
and of its derivatives of all orders, which were discussed in detail
in~\cite{CAoRFII}.

Another way to look at this issue is through the properties of the sets of
complex coefficients $c_{k}$ of the Taylor series. Given any set of
complex coefficients $c_{k}$, regardless of whether or not they follow
from a known inner analytic function, one can construct both a complex
power series $S(z)$ and the corresponding real coefficients $\alpha_{k}$
and $\beta_{k}$, using the relations in Equations~(\ref{EQTayCoe})
and~(\ref{EQPowSer}). In many cases the Fourier series generated by these
real coefficients will not converge, even if the complex power series
converges to an inner analytic function within the open unit disk.
However, if the complex power series is indeed convergent on that disk,
then one can discuss whether or not a real object can be defined on the
unit circle, through the $\rho\to 1_{(-)}$ limit from the open unit disk,
for example using the summation rule for Fourier series given in
Equation~(\ref{EQSummRule}).

If we examine that summation rule, it is apparent that it will work for
much more than just integrable real functions, which always have bounded
Fourier coefficients. For example, one may have unbounded Fourier
coefficients $\alpha_{k}$ and $\beta_{k}$, such as those of the $n^{\rm
  th}$ derivative of the delta ``function'', which diverge to infinity as
the power $k^{n}$ when $k\to\infty$, and still have a well-defined inner
analytic function, as was shown in detail in~\cite{CAoRFII}. In fact, one
can show that the summation rule can be used for all sets of Fourier
coefficients that do {\em not} diverge exponentially fast with $k$. In
order to develop this idea, let us first define a very general condition
on the sequences of complex coefficients that guarantees that the
corresponding power series are convergent within the open unit disk, and
thus converge to inner analytic functions.

\begin{definition}\Colon\label{Def01}
  Exponentially Bounded Coefficients
\end{definition}

\noindent
Given an arbitrary ordered set of complex coefficients $a_{k}$, for
$k\in\{0,1,2,3,\ldots,\infty\}$, if they satisfy the condition that

\noindent
\begin{equation}\label{EQExpLim1}
  \lim_{k\to\infty}
  |a_{k}|\e{-Ck}
  =
  0,
\end{equation}

\noindent
for all real $C>0$, then we say that the sequence of coefficients $a_{k}$
is {\em exponentially bounded}.

\vspace{2.6ex}

\noindent
What this means is that $a_{k}$ may or may not go to zero as $k\to\infty$,
may approach a non-zero complex number, and may even diverge to infinity
as $k\to\infty$, so long as it does not do so exponentially fast. This
includes therefore not only the sequences of complex Taylor coefficients
corresponding to all possible convergent Fourier series, but many
sequences that correspond to Fourier series that diverge almost
everywhere. Also, it not only includes the sequences of complex Taylor
coefficients corresponding to all possible integrable real functions, but
many sequences of coefficients that cannot be obtained at all from a real
function, such as those associated to the Dirac delta ``function'' and its
derivatives of arbitrarily high orders, as was shown in~\cite{CAoRFII}. We
see therefore that this is a very weak condition on the complex sequence
of coefficients $a_{k}$.

Before we proceed to the extension of the Fourier theory, let us establish
a preliminary result, which can be understood as a property of the
sequences of complex coefficients $a_{k}$ which satisfy the condition
stated in Definition~\ref{Def01}. We will show that the condition
expressed in Equation~(\ref{EQExpLim1}) implies an infinite collection of
other similar conditions involving the $k\to\infty$ limit, that express
modified bounds on these sequences of coefficients.

\begin{property}\Colon\label{Prop1}
  If the sequence of complex coefficients $a_{k}$ is exponentially
  bounded, then we also have that

  \noindent
  \begin{equation}\label{EQExpLim2}
    \lim_{k\to\infty}
    |a_{k}|k^{p}\e{-Ck}
    =
    0,
  \end{equation}

  \noindent
  for all real $C>0$ and for all real powers $p>0$.

\end{property}

\noindent
This is just a formalization of the well-known fact that the
negative-exponent real exponential function of $k$ goes to zero faster
than any positive real power of $k$ goes to infinity, as $k\to\infty$. In
order to prove this, we observe that for $k>0$ we may write the function
of $k$ on the left-hand side of Equation~(\ref{EQExpLim2}) as

\begin{equation}
  |a_{k}|k^{p}\e{-Ck}
  =
  |a_{k}|\e{\,p\ln(k)}\e{-Ck}.
\end{equation}

\noindent
Note that this is a positive real quantity. Recalling the properties of
the real logarithm function, we now observe that, given an arbitrary real
number $A>0$, there is always a sufficiently large finite value $k_{m}$ of
$k$ above which $\ln(k)<Ak$. Due to this we may write, for all $k>k_{m}$,

\begin{equation}
  |a_{k}|k^{p}\e{-Ck}
  <
  |a_{k}|\e{\,pAk}\e{-Ck},
\end{equation}

\noindent
since the exponential with a strictly positive real exponent is a
monotonically increasing function. If we now choose $A=C/(2p)$, which we
may do because this value is positive and not zero, we get that, for all
$k>k_{m}$,

\noindent
\begin{eqnarray}
  |a_{k}|k^{p}\e{-Ck}
  & < &
  |a_{k}|\e{Ck/2}\e{-Ck}
  \nonumber\\
  & = &
  |a_{k}|\e{-Ck/2}.
\end{eqnarray}

\noindent
According to our hypothesis about the coefficients $a_{k}$, the
$k\to\infty$ limit of the expression in the right-hand side is zero for
any strictly positive value of $C'=C/2$, so that taking the $k\to\infty$
limit we establish our preliminary result,

\begin{equation}
  \lim_{k\to\infty}
  |a_{k}|k^{p}\e{-Ck}
  =
  0,
\end{equation}

\noindent
for all real $C>0$ and all real $p>0$. Therefore, we have established this
property.

\vspace{2.6ex}

\noindent
Let us now show that the condition that the sequence of complex
coefficients $c_{k}$ in Equation~(\ref{EQTayCoe}) is exponentially bounded
is equivalent to the condition that the sequences of real coefficients
$\alpha_{k}$ and $\beta_{k}$ are both exponentially bounded. First, if we
assume that the sequences $\alpha_{k}$ and $\beta_{k}$ are both
exponentially bounded, and since from Equation~(\ref{EQTayCoe}) we have
that

\begin{equation}
  |c_{k}|
  =
  \sqrt{|\alpha_{k}|^{2}+|\beta_{k}|^{2}},
\end{equation}

\noindent
it follows at once that

\noindent
\begin{eqnarray}
  \lim_{k\to\infty}
  |c_{k}|\e{-Ck}
  & = &
  \lim_{k\to\infty}
  \sqrt{
    \left(|\alpha_{k}|\e{-Ck}\right)^{2}
    +
    \left(|\beta_{k}|\e{-Ck}\right)^{2}
  }
  \nonumber\\
  & = &
  \sqrt{
    \left(\lim_{k\to\infty}|\alpha_{k}|\e{-Ck}\right)^{2}
    +
    \left(\lim_{k\to\infty}|\beta_{k}|\e{-Ck}\right)^{2}
  }
  \nonumber\\
  & = &
  0,
\end{eqnarray}

\noindent
since both limits in the right-hand side are zero, thus establishing that
the sequence $c_{k}$ is exponentially bounded. Second, if we assume that
the sequence $c_{k}$ is exponentially bounded, and since from
Equation~(\ref{EQTayCoe}) we have that

\noindent
\begin{eqnarray}
  |c_{k}|
  & = &
  \sqrt{|\alpha_{k}|^{2}+|\beta_{k}|^{2}}
  \nonumber\\
  & \geq &
  |\alpha_{k}|
  \;\;\;\Rightarrow
  \nonumber\\
  |c_{k}|\e{-Ck}
  & \geq &
  |\alpha_{k}|\e{-Ck},
\end{eqnarray}

\noindent
taking the $k\to\infty$ limit and using the assumed property of the
sequence of coefficients $c_{k}$ it follows that

\noindent
\begin{eqnarray}
  \lim_{k\to\infty}
  |c_{k}|\e{-Ck}
  & \geq &
  \lim_{k\to\infty}
  |\alpha_{k}|\e{-Ck}
  \;\;\;\Rightarrow
  \nonumber\\
  0
  & \geq &
  \lim_{k\to\infty}
  |\alpha_{k}|\e{-Ck}
  \;\;\;\Rightarrow
  \nonumber\\
  \lim_{k\to\infty}
  |\alpha_{k}|\e{-Ck}
  & = &
  0,
\end{eqnarray}

\noindent
thus establishing that the sequence $\alpha_{k}$ is exponentially bounded.
Clearly, an identical argument can be made for the sequence $\beta_{k}$.
This establishes that the statement that the sequence of complex
coefficients $c_{k}$ is exponentially bounded is equivalent to the
statement that the sequences of real coefficients $\alpha_{k}$ and
$\beta_{k}$ are both exponentially bounded.

\vspace{2.6ex}

\noindent
Let us now prove the following theorem about the convergence of the power
series constructed out of a given arbitrary sequence of complex
coefficients $c_{k}$.

\begin{theorem}\Colon\label{Theo02}
  If the sequence of complex coefficients $c_{k}$, for
  $k\in\{0,1,2,3,\ldots,\infty\}$, is exponentially bounded, then the
  power series constructed from this sequence of coefficients converges
  within the open unit disk.
\end{theorem}

\noindent
Given the arbitrary sequence of complex coefficients $c_{k}$, we may
construct the complex power series in the complex $z$ plane, just as we
did in~\cite{CAoRFI},

\begin{equation}
  S(z)
  =
  \sum_{k=0}^{\infty}
  c_{k}z^{k}.
\end{equation}

\noindent
We will first show that, if the sequence of coefficients $c_{k}$ is
exponentially bounded, then this series is absolutely convergent inside
the open unit disk, which then implies that it is simply convergent there.

\begin{proof}\Colon
\end{proof}

\noindent
In order to prove that $S(z)$ is absolutely convergent, we consider the
real power series $\overline{S}(z)$ of the absolute values of the terms of
that series, which we write as

\noindent
\begin{eqnarray}
  \overline{S}(z)
  & = &
  \sum_{k=0}^{\infty}
  |c_{k}|\rho^{k}
  \nonumber\\
  & = &
  \sum_{k=0}^{\infty}
  |c_{k}|\e{k\ln(\rho)}.
\end{eqnarray}

\noindent
Since $\rho<1$ inside the open unit disk, the logarithm shown is strictly
negative, and we may put $\ln(\rho)=-C$ with real $C>0$. We can now see
that, according to our hypothesis about the coefficients $c_{k}$, the
terms of this series go to zero as $k\to\infty$,

\begin{equation}
  \overline{S}(z)
  =
  \sum_{k=0}^{\infty}
  |c_{k}|\e{-Ck},
\end{equation}

\noindent
since $C$ is real and strictly positive. In order to establish the
convergence of this real series, we write

\begin{equation}
  \overline{S}(z)
  =
  |c_{0}|
  +
  \sum_{k=1}^{\infty}
  \frac{k^{2}|c_{k}|\e{-Ck}}{k^{2}}.
\end{equation}

\noindent
According to the property expressed in Equation~(\ref{EQExpLim2}), with
$p=2$, the numerator shown above goes to zero as $k\to\infty$, and
therefore above a sufficiently large value $k_{m}$ of $k$ it is less that
one, so that we may write that

\noindent
\begin{eqnarray}
  \overline{S}(z)
  & = &
  \sum_{k=0}^{k_{m}}
  |c_{k}|\e{-Ck}
  +
  \sum_{k=k_{m}+1}^{\infty}
  \frac{k^{2}|c_{k}|\e{-Ck}}{k^{2}}
  \nonumber\\
  & < &
  \sum_{k=0}^{k_{m}}
  |c_{k}|\e{-Ck}
  +
  \sum_{k=k_{m}+1}^{\infty}
  \frac{1}{k^{2}}.
\end{eqnarray}

\noindent
The first term on the right-hand side is a finite sum and therefore is
finite, and the second term can be bounded from above by a convergent
asymptotic integral on $k$, so that we have

\noindent
\begin{eqnarray}
  \overline{S}(z)
  & < &
  \sum_{k=0}^{k_{m}}
  |c_{k}|\e{-Ck}
  +
  \int_{k_{m}}^{\infty}dk\,
  \frac{1}{k^{2}}
  \nonumber\\
  & = &
  \sum_{k=0}^{k_{m}}
  |c_{k}|\e{-Ck}
  +
  \frac{-1}{k}\at{k_{m}}{\infty}
  \nonumber\\
  & = &
  \sum_{k=0}^{k_{m}}
  |c_{k}|\e{-Ck}
  +
  \frac{1}{k_{m}}.
\end{eqnarray}

\noindent
This last expression is therefore a finite upper bound for all the partial
sums of the series $\overline{S}(z)$. It follows that $\overline{S}(z)$,
which is a real sum of positive terms, so that its partial sums form a
monotonically increasing real sequence which is now found to be bounded
from above, is therefore convergent. It then follows that $S(z)$ is
absolutely convergent and therefore convergent. Since this is valid for
all $\rho<1$, we may conclude that $S(z)$ converges on the open unit
disk. This completes the proof of Theorem~\ref{Theo02}.

\vspace{2.6ex}

\noindent
Since the series $S(z)$ considered above is a convergent power series
within the open unit disk, it converges to an analytic function $w(z)$ in
that domain, which is therefore an inner analytic function. We therefore
conclude that, if the sequence of complex coefficients $c_{k}$ in
Equation~(\ref{EQTayCoe}) is exponentially bounded, then it is the set of
Taylor coefficients of an inner analytic function. It now follows that, if
the corresponding Fourier coefficients $\alpha_{k}$ and $\beta_{k}$ are
both exponentially bounded, then the corresponding complex coefficients
$c_{k}$ are also exponentially bounded, and therefore the corresponding
Fourier series can be regulated by the use of the summation rule in
Equation~(\ref{EQSummRule}). Unless the Fourier coefficients go to zero as
$k\to\infty$, the Fourier series on the unit circle is sure to diverge
almost everywhere. One can then consider defining the corresponding real
object on the unit circle using the $\rho\to 1_{(-)}$ limit from the open
unit disk, for example through the use of the summation rule for the
Fourier series, given in Equation~(\ref{EQSummRule}).

In this way the Fourier theory of integrable real functions on the unit
circle can be extended to a much larger set of real objects, including for
example all the singular distributions discussed in~\cite{CAoRFII}, as
well as the examples of non-integrable real functions mentioned in that
paper. In fact, this extension of the Fourier theory includes a large
class of non-integrable real functions, as will be shown in the fourth
paper of this series. In this extended Fourier theory the real objects can
be considered as representable directly by their sequences of Fourier
coefficients, even when the corresponding Fourier series diverge. All
operations involving these divergent Fourier series can be mapped to
absolutely and uniformly convergent series and analytic operations within
the open unit disk, whose results are then taken to the unit circle
through the use of the $\rho\to 1_{(-)}$ limit. In many simple cases the
mere values of the real objects on the unit circle will be recovered in
this way, and in other more abstract cases global properties of the real
objects may be obtained in this way, such as in the case of the Dirac
delta ``function'' and its derivatives of all orders, as was discussed in
detail in~\cite{CAoRFII}.

\section{Conclusions and Outlook}\label{Sec07}

We have shown that the complex-analytic structure within the unit disk of
the complex plane established in a previous paper~\cite{CAoRFI}, which
leads to a close and deep relationship between integrable real functions
on the unit circle and inner analytic functions within the unit disk
centered at the origin of the complex plane, includes the whole structure
of the Fourier theory of integrable real functions. This fact leads to the
definition of a very general and powerful summation rule for Fourier
series, which allows one to still use and manipulate in a consistent way
divergent Fourier series, even when they are explicitly and strongly
divergent. The connection of the complex-analytic structure with the usual
Fourier theorems was exhibited.

The Fourier theory was then extended to include all the inner analytic
functions associated to singular Schwartz distributions, which were
discussed in detail in another previous paper~\cite{CAoRFII}, in which the
discussion of the complex-analytic structure was generalized to include
those singular distributions. In fact, the Fourier theory can be extended
to essentially the whole space of inner analytic functions. This includes
at least some non-integrable real functions, as was pointed out
in~\cite{CAoRFII}. The generalization to a much wider class of
non-integrable real functions will be tackled in a future paper.

As part of this process of extension, we introduced the concept of an
exponentially bounded sequence of complex coefficients $c_{k}$, and proved
that any such sequence is the set of Taylor coefficients of some inner
analytic function. As interesting open question is whether or not the
reverse of this statement is true, that is, whether or not the criterion
that the sequences of complex coefficients of the power series be
exponentially bounded includes all possible inner analytic functions. At
this time this seems rather unlikely, and in that case the problem poses
itself of what more general condition on the coefficients could cover the
whole space of inner analytic functions.

We believe that the results presented here establish a new perspective for
the study of the Fourier theory of real functions and related objects. It
provides a simple and complete account of all the mathematical structures
involved, as well as of all the main results of that theory, including in
particular a simple and solid proof of the completeness of the basis. Due
to this, it might also constitute a simpler and more efficient way to
teach the subject.

\section*{Acknowledgments}

The author would like to thank his friend and colleague Prof. Carlos
Eugênio Imbassay Carneiro, to whom he is deeply indebted for all his
interest and help, as well as his careful reading of the manuscript and
helpful criticism regarding this work.

\appendix

\section{Appendix: Scalar Product for Inner Analytic
  Functions}\label{App01}

Given two inner analytic functions $w_{1}(z)$ and $w_{2}(z)$, we consider
the complex contour integral over the circle $C_{0}$ of radius $\rho_{0}$,
with $0<\rho_{0}<1$, given by

\begin{equation}
  (w_{1}|w_{2})
  =
  \frac{1}{2\pi\ii}
  \oint_{C_{0}}dz\,
  \frac{1}{z}\,
  w_{1}^{*}(z)w_{2}(z).
\end{equation}

\noindent
Since the integrand in this expression is {\em not} analytic, the integral
depends on the circuit, and therefore on $\rho_{0}$. Therefore, what we
have here is in fact a one-parameter family of integrals. We will show
that for each value of $\rho_{0}$ this integral defines a scalar product
within the space of inner analytic functions, which induces in that space
a positive-definite norm. If we write the integral in terms of the
integration variable $\theta$, with constant $\rho_{0}$, we get for this
scalar product

\begin{equation}
  (w_{1}|w_{2})
  =
  \frac{1}{2\pi}
  \int_{-\pi}^{\pi}d\theta\,
  w_{1}^{*}(\rho_{0},\theta)w_{2}(\rho_{0},\theta).
\end{equation}

\noindent
If we now make both $w_{1}(z)$ and $w_{2}(z)$ equal to
$w(z)=u(\rho,\theta)+\ii v(\rho,\theta)$, we get

\noindent
\begin{eqnarray}
  (w|w)
  & = &
  \|w\|^{2}
  \nonumber\\
  & = &
  \frac{1}{2\pi}
  \int_{-\pi}^{\pi}d\theta\,
  |w(\rho_{0},\theta)|^{2}
  \nonumber\\
  & = &
  \frac{1}{2\pi}
  \int_{-\pi}^{\pi}d\theta\,
  \left[
    u^{2}(\rho_{0},\theta)
    +
    v^{2}(\rho_{0},\theta)
  \right]
  \nonumber\\
  & \geq &
  0,
\end{eqnarray}

\noindent
which is a manifestly real and positive quantity, that is zero if and only
if $w(\rho_{0},\theta)=0$ for all $\theta$, which in turn is equivalent to
$w(\rho,\theta)=0$ for all $\theta$ and all $\rho$ within the open unit
disk, because all zeros of an analytic function must be isolated, unless
it is the identically zero function. Therefore, for each value of the
parameter $\rho_{0}$ the real quantity $\|w\|$ is a positive-definite norm
on the space of all inner analytic functions which, as was observed
in~\cite{CAoRFI}, forms a vector space over the field of complex numbers.
That vector space is thus seen to constitute a complex Hilbert space, with
this scalar product and the associated positive-definite norm.

We can also see from the equation above that the scalar product and the
norm reduce naturally to the corresponding definitions for the real
functions $u(1,\theta)$ and $v(1,\theta)$ on the unit circle, when we take
the $\rho_{0}\to 1_{(-)}$ limit, thus establishing a close correspondence
between these two identical real Hilbert spaces on the unit circle and the
complex Hilbert space on the unit disk. In addition to this, for any value
of $\rho_{0}$ within the open interval $(0,1)$ we also have a pair of
identical real Hilbert spaces with the real functions $u(\rho_{0},\theta)$
and $v(\rho_{0},\theta)$ on the circle of radius $\rho_{0}$.

We may now show that the Taylor basis of functions around the origin,
which is complete to generate the whole space of inner analytic functions,
and which consists of the set of non-negative powers

\begin{equation}
  \left\{
    \rule{0em}{2.5ex}z^{k},
    k\in\{0,1,2,3,\ldots,\infty\}
  \right\},
\end{equation}

\noindent
is in fact an orthogonal basis according to this definition of the scalar
product. If we make $w_{1}(z)=w_{k_{1}}(z)=z^{k_{1}}$ and
$w_{2}(z)=w_{k_{2}}(z)=z^{k_{2}}$, we get

\noindent
\begin{eqnarray}
  (w_{k_{1}}|w_{k_{2}})
  & = &
  \frac{1}{2\pi}
  \int_{-\pi}^{\pi}d\theta\,
  \left(z^{k_{1}}\right)^{*}z^{k_{2}}
  \nonumber\\
  & = &
  \rho_{0}^{k_{1}+k_{2}}\,
  \frac{1}{2\pi}
  \int_{-\pi}^{\pi}d\theta\,
  \e{-\iii k_{1}\theta}\e{\iii k_{2}\theta}.
\end{eqnarray}

\noindent
Using now the first result shown in Equation~(\ref{EQComOrt}) we obtain
the orthogonality relation for the Taylor basis,

\begin{equation}
  (w_{k_{1}}|w_{k_{2}})
  =
  \rho_{0}^{k_{1}+k_{2}}\,
  \delta_{k_{1},k_{2}}.
\end{equation}

\noindent
Since the integer powers are analytic on the whole complex plane, there is
no obstruction to taking the $\rho_{0}\to 1_{(-)}$ limit, and thus we see
that in this case the Taylor basis is not only orthogonal, but also
normalized,

\begin{equation}
  (w_{k_{1}}|w_{k_{2}})
  =
  \delta_{k_{1},k_{2}},
\end{equation}

\noindent
with $\|w_{k}\|=1$ for all $k$, where the scalar product is now defined on
the unit circle. If we write the inner analytic functions in terms of
their Taylor series around the origin,

\noindent
\begin{eqnarray}
  w_{1}(z)
  & = &
  \sum_{k=0}^{\infty}
  c_{1,k}z^{k},
  \nonumber\\
  w_{2}(z)
  & = &
  \sum_{k=0}^{\infty}
  c_{2,k}z^{k},
\end{eqnarray}

\noindent
we obtain for the scalar product, since we may always integrate convergent
power series term-by-term,

\noindent
\begin{eqnarray}
  (w_{1}|w_{2})
  & = &
  \frac{1}{2\pi\ii}
  \oint_{C_{0}}dz\,
  \frac{1}{z}
  \sum_{k_{1}=0}^{\infty}
  \sum_{k_{2}=0}^{\infty}
  c_{1,k_{1}}^{*}c_{2,k_{2}}
  \left(z^{k_{1}}\right)^{*}z^{k_{2}}
  \nonumber\\
  & = &
  \sum_{k_{1}=0}^{\infty}
  \sum_{k_{2}=0}^{\infty}
  c_{1,k_{1}}^{*}c_{2,k_{2}}\,
  \frac{1}{2\pi\ii}
  \oint_{C_{0}}dz\,
  \frac{1}{z}
  \left(z^{k_{1}}\right)^{*}z^{k_{2}}
  \nonumber\\
  & = &
  \sum_{k_{1}=0}^{\infty}
  \sum_{k_{2}=0}^{\infty}
  c_{1,k_{1}}^{*}c_{2,k_{2}}
  (w_{k_{1}}|w_{k_{2}})
  \nonumber\\
  & = &
  \sum_{k_{1}=0}^{\infty}
  \sum_{k_{2}=0}^{\infty}
  c_{1,k_{1}}^{*}c_{2,k_{2}}\,
  \rho_{0}^{k_{1}+k_{2}}\,
  \delta_{k_{1},k_{2}}
  \nonumber\\
  & = &
  \sum_{k=0}^{\infty}
  \rho_{0}^{2k}
  c_{1,k}^{*}c_{2,k},
\end{eqnarray}

\noindent
where we identified the scalar product $(w_{k_{1}}|w_{k_{2}})$ and then
used the orthogonality relations of the Taylor basis. So long as
$\rho_{0}<1$, and so long as $c_{1,k}$ and $c_{2,k}$ are exponentially
bounded, this series converges exponentially fast. We may also write the
corresponding expression for the norm, if we make $c_{1,k}=c_{2,k}=c_{k}$
and $w_{1}(z)=w_{2}(z)=w(z)$,

\noindent
\begin{eqnarray}
  \|w\|^{2}
  & = &
  (w|w)
  \nonumber\\
  & = &
  \sum_{k=0}^{\infty}
  \rho_{0}^{2k}
  |c_{k}|^{2},
\end{eqnarray}

\noindent
with the same conditions for the convergence of the series. In all this
structure, if we take the $\rho_{0}\to 1_{(-)}$ limit, the scalar product
and the norm may in general diverge, unlike what happens in the case of
the elements of the Taylor basis. However, so long as $\rho_{0}<1$ all the
inner analytic functions have finite norms and finite scalar products with
one another. In some cases, it may be possible to determine the values of
these quantities on the unit circle using the $\rho_{0}\to 1_{(-)}$ limit,
even if the corresponding series expressions written directly on the unit
circle diverge.

Perhaps the best way to characterize this structure is as a one-parameter
family of pairs of identical real Hilbert spaces, one associated to the
real parts and another associated to the imaginary parts of the inner
analytic functions, where the parameter is the radius $\rho_{0}$ of each
circle within the unit disk, which are connected to each other by a
process of analytic continuation. For each value of $\rho_{0}$ within the
open interval $(0,1)$ there is a one-to-one mapping between the inner
analytic functions on the open unit disk and the real functions obtained
as the real parts of these inner analytic function restricted to the
circle of radius $\rho_{0}$. This one-to-one mapping preserves the scalar
product and the norm, as they are defined within each space. This fact is
still true even in the $\rho_{0}\to 1_{(-)}$ limit, although in that case
not every real object at the unit circle, resulting from the limit, is a
normal real function, and although in many cases the norms and scalar
products may diverge in the limit.

Note that the integral defining the scalar product of the inner analytic
functions is a one-dimensional integral over the circle of radius
$\rho_{0}$, despite the fact that each complex inner analytic function
consist of a pair of real functions of two variables. However, this is a
natural characteristic of the scalar product in this context, since it is
a well-known fact that an analytic function is completely determined on a
two-dimensional region of the complex plane by its values only at a
one-dimensional boundary of that region. In this way, although only a
one-dimensional restriction of the inner analytic function is explicitly
taken into account in the integral over the circle of radius $\rho_{0}$
that defines the scalar product, that restriction still includes
implicitly the whole structure of the inner analytic function within the
corresponding disk of radius $\rho_{0}$. Therefore, it is perhaps arguable
that the most natural definition of the scalar product is that associated
to the choice $\rho_{0}=1$, despite the convergence issues that this
choice may involve.

\bibliography{allrefs_en}\bibliographystyle{ieeetr}

\end{document}

%% file: packages.tex
\usepackage[latin1]{inputenc}
\usepackage{epsfig}
\usepackage{amsfonts}
\usepackage{cite}
\usepackage{color}
\definecolor{lgrey}{gray}{0.99}

%% file: math-defs.tex
\newcommand{\ii}{\mbox{\boldmath$\imath$}}
\newcommand{\iii}{\mbox{\boldmath\scriptsize$\imath$}}

\newcommand{\at}[2]{\left.\rule{0em}{3ex}\right[_{\,#1}^{\,#2}}

\newcommand{\e}[1]{\,{\rm e}^{#1}}

\newcommand{\ldot}{\mbox{\Large$\cdot$}\!}

\newtheorem{definition}{Definition}
\newtheorem{property}{Property}[definition]

\newtheorem{theorem}{Theorem}

\newtheorem{proof}{Proof}[theorem]
\newcommand{\Colon}{{\hspace{-0.4em}\bf:}}